\documentclass[12pt,leqno]{amsart}
\usepackage{amssymb}
\usepackage{amsmath,amssymb,color}
\numberwithin{equation}{section}

\newcommand{\ep}{\varepsilon}
\newcommand{\la}{\lambda}
\newcommand{\va}{\varphi}
\newcommand{\ppp}{\partial}

\newcommand{\www}{\widetilde}

\newcommand{\R}{\mathbb{R}}

\newcommand{\N}{\mathbb{N}}

\newcommand{\ooo}{\overline}
\newcommand{\OOO}{\Omega}

\newcommand{\sumt}{\sum_{\vert \beta\vert\le 2}}
\newcommand{\sums}{\sum_{\vert \beta\vert\le 3}}
\newcommand{\dddO}{\Vert d\Vert_{C(\ooo{\OOO})}}
\newcommand{\dddOO}{\Vert d\Vert_{C(\ooo{\OOO_1})}}
\newcommand{\weight}{e^{2s\alpha}}

\newcommand{\hhalf}{\frac{1}{2}}

\allowdisplaybreaks

\title
[]
{
Inverse source problem and the continuation for a fourth-order 
parabolic equation in general dimensions
}

\author{
$^1$ O.~Yu.~Imanuvilov and \, $^{2,3,4}$ M.~Yamamoto }

\thanks{
$^1$ Department of Mathematics, Colorado State
University, 101 Weber Building, Fort Collins, CO 80523-1874, U.S.A.
e-mail: {\tt oleg@math.colostate.edu}\\
$^2$ Graduate School of Mathematical Sciences, The University
of Tokyo, Komaba, Meguro, Tokyo 153-8914, Japan \\
$^3$ Honorary Member of Academy of Romanian Scientists, 
Splaiul Independentei Street, no 54,
050094 Bucharest Romania \\
$^4$ Peoples' Friendship University of Russia 
(RUDN University) 6 Miklukho-Maklaya St, Moscow, 117198, Russian Federation
e-mail: {\tt myama@ms.u-tokyo.ac.jp}
}

\date{}
\begin{document}
\maketitle

\baselineskip 18pt

\begin{abstract}
In this article, for a fourth-order parabolic equation which is 
closely related for example to the Cahn-Hilliard equation,
we study an inverse source problem by interior data 
and the continuation of solution 
from lateral Cauchy data.  Our method relies on a Carleman estimate and 
proves conditional stability for both problems.
\\
{\bf Key words.}  
Fourh-order parabolic equation,
inverse source problem, continuation, Carleman estimate, stability
\\
{\bf AMS subject classifications.}
35R30, 35R25
\end{abstract}

\section{Introduction and main results}
In this article, let $\OOO \subset \R^n$, $n\ge 2$ be an open bounded 
domain with smooth boundary $\ppp\OOO$ and let $T>0$.
By $\nu(x)$ we denote the unit outward normal vector to $\ppp\OOO$ at
$x=(x_1,\dots,x_n)$.  We set $\ppp_{\nu}$ as pointwised normal derivative and also
as the trace in Sobolev spaces.
For $\beta := (\beta_1, ..., \beta_n) \in (\N \cup \{0\})^n$, we put
$$
\ppp_x^{\beta} = \ppp_1^{\beta_1}\cdots \ppp_n^{\beta_n},\quad 
\vert \beta\vert = \beta_1 + \cdots + \beta_n, \quad 
\ppp_k = \frac{\ppp}{\ppp x_k}, \quad k=1,2,..., n.
$$

Mainly we consider
$$
P(t,x,D)y :=
\ppp_ty + \Delta^2y(x,t) + \sumt p_{\beta}(x)\ppp_x^{\beta}y
= F(x,t), \quad x\in \OOO,\, 0<t<T
                                                           \eqno{(1.1)}
$$
with $p_{\beta} \in L^{\infty}(\OOO)$.

Equation (1.1) is related to the Cahn-Hilliard equation and see
Guererro and Kassab \cite{GK} as for other applications, which 
studies the null controllability.

In this article, we consider two   problems.

{\bf (I): Inverse source problem}.\\
{\it Let $\theta \in (0,T)$ and an open  subdomain $\omega \subset \OOO$ be 
arbitrarily given.  For
$$
\ppp_ty + \Delta^2y(x,t) + \sumt p_{\beta}(x)\ppp_x^{\beta}y
= R(x,t)f(x), \quad x\in \OOO,\, 0<t<T
                                                           \eqno{(1.2)}
$$
with known $R$ and 
$$
y = \Delta y = 0 \quad \mbox{on $\ppp\OOO \times (0,T)$,}
                                                         \eqno{(1.3)}
$$
we need to determine $f=f(x)$, $x\in \OOO$ by data
$$
y\vert_{\omega\times (0,T)} \quad \mbox{and} \quad y(\cdot,\theta)
\quad \mbox{in $\OOO$}.
$$ }

{\bf (II): Continuation of solution}.\\
Let $\Gamma \subset \ppp\OOO$ be an arbitrarily given subboundary and let
$$
\ppp_tu + \Delta^2u(x,t) + \sumt p_{\beta}(x)\ppp_x^{\beta}u
= 0, \quad x\in \OOO,\, 0<t<T.
                                                           \eqno{(1.4)}
$$
Let $\OOO_0 \subset \OOO$ be a open subdomain such that 
$\ooo{\OOO_0} \subset \OOO \cup \Gamma$ and let $\ep>0$ be arbitrary.
Then determine $u\vert_{\OOO_0 \times (\ep, T-\ep)}$ by 
$\ppp_{\nu}^ju\vert_{\Gamma \times (0,T)}$, $0\le j \le 3$.
\\

Now we state the main result for Problem (I):\\
{\bf Theorem 1.}\\
{\it
Let $y \in H^1(0,T;H^4(\OOO)) \cap H^2(0,T;L^2(\OOO))$ satisfy (1.2) and (1.3),
and let there exist a constant $r_0 > 0$ such that 
$$
R, \, \ppp_tR \in L^2(0,T;L^{\infty}(\OOO)), \quad 
\vert R(x,\theta)\vert \ge r_0, \quad x\in \ooo{\OOO}.
                                                            \eqno{(1.5)}
$$
We arbitrarily fix $\theta\in (0,T)$ and chose  $t_1 > 0$ such that 
$0 < \theta-t_1<\theta + t_1 < T$.
Then there exists a constant $C>0$ such that 
$$
\Vert f\Vert_{L^2(\OOO)} \le C(\Vert y\Vert_{H^1(\theta-t_1,\theta+t_1;
L^2(\omega))}
+ \Vert y(\cdot,\theta)\Vert_{H^4(\OOO)}).
$$
}
\\

This theorem asserts the global Lipschitz stability in determining 
a spatially varying factor $f(x)$ of the source term $R(x,t)f(x)$ and we use
extra data $y\vert_{\omega \times (\theta-t_1,\theta+t_1)}$ and 
$y(\cdot,\theta)$, provided that we are given boundary data 
$y$ and $\Delta y$ on the whole lateral boundary $\ppp\OOO\times 
(\theta-t_1,\theta+t_1)$ and $\theta \in (0,T)$.

This type of the global Lipschitz stability was established for the
parabolic inverse source problem (Imanuvilov and Yamamoto \cite{IY1}).
We remark that we do not know whether the global Lipschitz stability 
holds for the case of $\theta=0$, and we may conjecture negatively.

As for a similar inverse problem for a one-dimensional 
fourth-order parabolic equation, we can refer to 
Baudouin, Cerpa, Cr\'epeau and Mercado \cite{BCCM}.
\\

Next we show the second main result for Problem (II).
\\
{\bf Theorem 2.}\\
{\it
We assume that $u \in L^2(0,T;H^4(\OOO)) \cap H^1(0,T;L^2(\OOO))$ 
satisfies (1.4).
Let $\Gamma \subset \ppp\OOO$ be any given subboundary and 
let $\OOO_0 \subset \OOO$ be an open subdomain such that 
$\ooo{\OOO_0} \subset \OOO \cup \Gamma$ and let $\ep\in (0,T), \mu>0$ be 
arbitrary.
Then there exist constants $C>0$ and $\kappa \in (0,1)$ such that
\begin{align*}
&\Vert u\Vert_{H^1(\ep,T-\ep;L^2(\OOO_0))}
+ \Vert u\Vert_{L^2(\ep,T-\ep;H^2(\OOO_0))}
+ \Vert \nabla (\Delta u)\Vert_{L^2(\ep,T-\ep;H^2(\OOO_0))}
+ \Vert \Delta^2 u\Vert_{L^2(\ep,T-\ep;H^2(\OOO_0))}
\\
\le& C_0(\Vert u\Vert_{L^2(0,T;H^3(\OOO))})
(D + D^{\kappa}),
\end{align*}
where we set 
$$
D:= \sum_{j=0}^3 \Vert \ppp_{\nu}^ju\Vert
_{L^2(0,T;H^{\frac{7}{2}-j}(\Gamma))}
+ \Vert u\Vert_{H^1(0,T;H^{\mu}(\Gamma))}
$$
and the constant $C_0(\eta)>0$ depends on $\eta>0$.
}
\\

Theorem 2 provides a conditional stability estimate of H\"older type
in continuing a solution: if data $(u,\ppp_{\nu}u,\ppp_{\nu}^2u,\ppp_{\nu}^3u)$ on 
$\Gamma \times (0,T)$ are small, 
then we can estimate $u$ in a subdomain of $\OOO\times (0,T)$, 
provided that $\Vert u\Vert_{L^2(0,T;H^3(\OOO))}$ is bounded.
We note that we need not assume any data outside of
$\Gamma \times (0,T)$.

We cannot choose $\OOO_0 = \OOO$ and $\ep=0$ keeping the 
H\"older stability.  However, since $\OOO_0 \subset \OOO$ and 
$\ep>0$ are arbitrary in Theorem 2, we can easily derive the 
uniqueness.
\\
{\bf Corollary.}\\
{\it
If $u \in L^2(0,T;H^4(\OOO)) \cap H^1(0,T;L^2(\OOO))$ satisfy (1.4)
and 
$$
\ppp_{\nu}^j u = 0 \quad \mbox{on}\quad  \Gamma \times (0,T),\quad 
\forall j\in \{0,1,2,3\},
$$
then $u=0$ in $\OOO \times (0,T)$.
}
\\

Our proofs are based on Carleman estimate, see Bukhgeim and Klibanov
\cite{BK}, Klibanov \cite{Kl2} as pioneering works.
We refer for example to Bellassoued and Yamamoto \cite{BY},
Imanuvilov and Yamamoto \cite{IY2}, Klibanov and Timonov \cite{KT}.
 
The article is composed of five sections.  In Section 2, we establish  
a Carleman estimate.  Sections 3 and 4 are devoted to the proofs of 
Theorems 1 and 2 respectively.
In Section 5, we give concluding remarks.

\section{Key Carleman estimate.}

We present a Carleman estimate established in Guerrero and Kassab 
\cite{GK}.  Let $\omega_0 \subset \omega$ be a non-empty domain such 
that $\ooo{\omega_0} \subset \omega$.  First we introduce 
a function $d \in C^4(\ooo{\OOO})$ such that 
$$
d\vert_{\ppp\OOO} = 0, \quad \vert \nabla d(x)\vert > 0 \quad 
\mbox{on $\ooo{\OOO \setminus \omega_0}$}.     \eqno{(2.1)}
$$
The existence of function $d$ is proved in Imanuvilov \cite{Ima}.

Let  $\la>0$ be a positive parameter, we chose $\tau\in (0,\frac T2)$ and  $t_0>0$ such that $\tau<t_0<T-\tau.$ Next we set
\begin{align*}
& \alpha(x,t) = \alpha_{\tau,t_0}(x,t) 
= \frac{e^{\la d(x)} - e^{2\la\dddO}}{\sqrt{(t-(t_0-\tau))
(t_0+\tau-t)}},\\
& \va(x,t) = \va_{\tau,t_0}(x,t) 
= \frac{e^{\la d(x)}}{\sqrt{(t-(t_0-\tau))(t_0+\tau-t)}},
\quad (x,t) \in \OOO \times (t_0-\tau, \, t_0+\tau).
\end{align*}

Then
\\
{\bf Proposition 1 (Carleman estimate \cite{GK}).}\\
{\it 
Let $\la > 0$ be chosen sufficiency large.  Then there exist 
constants $s_0>0$ and $C>0$ such that 
\begin{align*}
& \int_{\OOO \times (t_0-\tau,t_0+\tau)}
\biggl( s^6\va^6\vert y\vert^2 + s^4\va^4\vert \nabla y\vert^2
+ s^2\va^2\vert \nabla (\nabla y)\vert^2
+ s\va\vert \nabla (\Delta y)\vert^2\\
+& s^{-1}\va^{-1}(\vert \ppp_ty\vert^2 + \vert \nabla^2y\vert^2)
\biggr) \weight dxdt
\end{align*}
$$
\le C\left( \int_{\OOO \times (t_0-\tau,t_0+\tau)} \vert P(t,x,D)y\vert^2\weight 
dxdt
+ \Vert y\Vert^2_{L^2(\omega \times (t_0-\tau,t_0+\tau))}\right)
                                                        \eqno{(2.2)}
$$
for all $s \ge s_0$ and  function $y\in L^2(t_0-\tau, t_0+\tau; H^4(\OOO))
\cap H^1(t_0-\tau, t_0+\tau; L^2(\OOO))$ satisfying (1.3).
}
\\

The proof is given in \cite{GK}.
Here and henceforth $C$ denotes generic, strictly positive  constants which are independent of parameter
$s.$  The constants $C>0$ and $s_0>0$ are 
independent of choices of the coefficients $p_{\beta}$ of 
(1.1), and depends on a bound $M_0$:
$\Vert p_{\beta}\Vert_{L^{\infty}(\OOO)} \le M_0$ for $\vert \beta \vert 
\le 2$.
Moreover, the constants  $C$ and $s_0$ are dependent  on $\tau$, but 
independent of 
$t_0 \in (\tau, T-\tau)$, because (1.1) and the Carleman estimate are 
invariant by the translation $t \longrightarrow t+t_1$ with any 
constant $t_1$ as long as $(t_0-\tau+t_1,\, t_0+\tau+t_1) \subset 
(0,T)$.
We can include other derivatives in $x$ on the left-hand side of (2.2)
but we omit.

This type of Carleman estimate is global in the sense that 
it holds over $\OOO$ in the $x$-direction and it was firstly proved in 
Imanuvilov \cite{Ima} for parabolic equations.
In \cite{GK}, the weight function:
$$
\www{\alpha}(x,t)  
= \frac{e^{2\la\dddO}(e^{\la d(x)} - e^{2\la\dddO})}
{\sqrt{(t-(t_0-\tau))(t_0+\tau-t)}}
$$
is used and different from ours, but replacing 
$$
\www{s}:= se^{2\la\dddO},
$$
we derive (2.2).  Moreover, in \cite{GK}, the term of 
$y\vert_{\omega\times (t_0-\tau,t_0+\tau)}$ is given by the integral
$$
\int_{\omega\times (t_0-\tau, t_0+\tau)}
s^7\va^7 \vert y\vert^2 \weight dxdt,
$$
but here we substitute the following:
setting 
$$
h(t) := \frac{1}{\sqrt{(t-(t_0-\tau))(t_0+\tau-t)}},
$$
we have
\begin{align*}
& s^7\va^7 e^{2s\alpha(x,t)}
= s^7h(t)^7e^{7\la d(x)}
\exp(-2s(e^{2\la\dddO} - e^{\la d(x)})h(t))\\
\le& C(sh(t))^7e^{-2C(sh(t))}
\le C\sup_{\eta\ge 0} \eta^7e^{-2C\eta} < \infty.
\end{align*}

\section{Proof of Theorem 1.}

We follow the arguments in \cite{IY1}, which considers a similar inverse 
source problem for a second-order parabolic equation on the basis of 
the relevant Carleman estimate.

Setting $a = y(\cdot,\theta)$ and $z = \ppp_ty$, we have
$$\left\{ \begin{array}{rl}
& \ppp_tz + \Delta^2 z + \sumt p_{\beta}(x)\ppp_x^{\beta} z
= (\ppp_tR)(x,t)f(x), \quad x\in \OOO, \, 0<t<T, \\
& z = \Delta z = 0 \quad \mbox{on $\ppp\OOO \times (0,T)$}
\end{array}\right.
                                              \eqno{(3.1)}
$$
and
$$
z(x,\theta) = -\Delta^2a - \sumt p_{\beta}\ppp_x^{\beta}a
+ R(x,\theta)f(x), \quad x \in \OOO.              \eqno{(3.2)}
$$
In Proposition 1, we set 
$$
t_0=\theta, \quad \tau = t_1,
$$
and write $\alpha = \alpha_{t_1,\theta}$ and
$\va = \va_{t_1,\theta}$.
Thus, applying Carleman estimate (2.2) to system (3.1), we obtain
$$
\int_{\OOO\times (\theta-t_1,\theta+t_1)} 
\left( s^6\va^6\vert z\vert^2 + \frac{1}{s\va}\vert \ppp_tz\vert^2
\right) \weight dxdt 
\le C\int_{\OOO\times (\theta-t_1,\theta+t_1)} 
\vert f\vert^2 \weight dxdt + CD_0^2           \eqno{(3.3)}
$$
for all large $s>0$.
Here we set
$$
D_0 := \Vert y\Vert_{H^1(\theta-t_1,\theta+t_1;L^2(\omega))}.
$$
On the other hand, since $\alpha(x,\theta-t_1) = -\infty$, we see
\begin{align*}
& \int_{\OOO} \vert z(x,\theta)\vert^2 e^{2s\alpha(x,\theta)} dx 
= \int^{\theta}_{\theta-t_1} \ppp_t
\left( \int_{\OOO} \vert z(x,t)\vert^2 \weight dx \right) dt\\
=& \int^{\theta}_{\theta-t_1}\int_{\OOO}
(2z\ppp_tz + 2s(\ppp_t\alpha)\vert z\vert^2) \weight dxdt.
\end{align*}
Provided that the parameter $\lambda$ is sufficiently large, one can directly 
verify
$$
\vert \ppp_t\alpha(x,t)\vert \le C\va^2(x,t), \quad 
(x,t) \in \OOO\times (\theta-t_1,\theta+t_1),
$$
and so 
$$
\int_{\OOO} \vert z(x,\theta)\vert^2 e^{2s\alpha(x,\theta)} dx 
\le C\int_{\Omega \times (\theta-t_1, \theta+t_1)}
(\vert z\vert \vert \ppp_tz\vert + s\va^2 \vert z\vert^2) \weight dxdt.
                                                       \eqno{(3.4)}
$$
Since 
$$
\vert z\vert \vert \ppp_tz\vert 
= s^{\hhalf}\va^{\hhalf}\vert z\vert s^{-\hhalf}\va^{-\hhalf}
\vert \ppp_tz\vert
\le \hhalf(s\va \vert z\vert^2 + s^{-1}\va^{-1}\vert \ppp_tz\vert^2)
$$
and $s\va^2\vert z\vert^2 \le Cs^6\va^6\vert z\vert^2$, we can apply
(3.3) to the right-hand side of inequality (3.4), so that 
$$
\int_{\OOO} \vert z(x,\theta)\vert^2 e^{2s\alpha(x,\theta)} dx 
\le C\int_{\OOO\times (\theta-t_1,\theta+t_1)} \vert f(x)\vert^2
\weight dxdt + CD_0^2
$$
for all large $s>0$.
By (3.2) and (1.5), we obtain
$$
\int_{\OOO} \vert f(x)\vert^2 e^{2s\alpha(x,\theta)} dx
\le C\Vert a\Vert^2_{H^4(\OOO)} 
+ C\int_{\OOO} \vert z(x,\theta)\vert^2 e^{2s\alpha(x,\theta)} dx 
$$
$$
\le C\int_{\OOO\times (\theta-t_1,\theta+t_1)} \vert f(x)\vert^2
\weight dxdt + C\www{D}^2                   \eqno{(3.5)}
$$
for all large $s>0$.  Here and henceforth we set 
$\www{D} = \Vert a\Vert_{H^4(\OOO)} + D_0$ and
$$
h(t) := \frac{1}{\sqrt{(t-(\theta-t_1))(\theta+t_1-t)}}.
$$
Now we can directly verify
$$
\alpha(x,t) - \alpha(x,\theta)
= (h(t)-h(\theta)) (e^{\la d(x)} - e^{2\la\dddO})
$$
and $h(t) > h(\theta)$ for $t \ne \theta$.
Hence
$$
\alpha(x,t) - \alpha(x,\theta)
\le -(h(t)-h(\theta)) (e^{2\la \dddO} - e^{\la\inf_{x\in\Omega}d(x)}), \quad
\theta-t_1 \le t \le \theta + t_1,
$$
and so 
\begin{align*}
& \int_{\OOO \times (\theta-t_1,\theta_1+t_1)} 
\vert f(x)\vert^2 e^{2s\alpha(x,t)} dxdt
= \int_{\OOO} \vert f(x)\vert^2 e^{2s\alpha(x,\theta)} 
\left( \int^{\theta+t_1}_{\theta-t_1} 
e^{2s(\alpha(x,t) - \alpha(x,\theta))} dt\right)  dx\\
\le & \int_{\OOO} \vert f(x)\vert^2 e^{2s\alpha(x,\theta)} 
\left( \int^{\theta+t_1}_{\theta-t_1} 
e^{-C_0s(h(t) - h(\theta))} dt\right)  dx.
\end{align*}
The Lebesgue convergence theorem implies 
$$
\int^{\theta+t_1}_{\theta-t_1} e^{-C_0s(h(t) - h(\theta))} dt
= o(1)\quad \mbox{as}\quad s\rightarrow +\infty.
$$
 Therefore 
$$
\int_{\OOO \times (\theta-t_1,\theta_1+t_1)} 
\vert f(x)\vert^2 e^{2s\alpha(x,t)} dxdt
= o(1) \int_{\OOO} \vert f(x)\vert^2 e^{2s\alpha(x,\theta)} dx.
$$
Substituting this into (3.5) and absorbing the term on the right-hand side 
into the left-hand side, we complete the proof of Theorem 1. $\square$

\section{Proof of Theorem 2.}

The proof is similar to the proof of Theorem 5.1 in Yamamoto \cite{Ya}.

{\bf First Step.}

We change the boundary data to zero and make some extensions to a wider 
spatial domain.

By the Sobolev extension theorem, we can find $\www{u} 
\in L^2(0,T;H^4(\OOO)) \cap H^1(0,T;L^2(\OOO))$ such that 
$$
\ppp_{\nu}^j\www{u} = \ppp_{\nu}^ju \quad \mbox{on}\quad  \Gamma \times (0,T)\quad \forall j\in \{0,\dots ,3\}
$$
and
$$
\Vert \www{u}\Vert_{L^2(0,T;H^4(\OOO))}
+ \Vert \ppp_t\www{u}\Vert_{L^2(0,T;L^2(\OOO))}
$$
$$
\le C\left( \sum_{j=0}^3 \Vert \ppp_{\nu}^ju\Vert
_{L^2(0,T;H^{\frac{7}{2}-j}(\Gamma))}
+ \Vert u\Vert_{H^1(0,T;H^{\mu}(\Gamma))}\right) 
=: CD.                                    \eqno{(4.1)}
$$
We set 
$$
v:= u - \www{u}.
$$
Then
$$
\ppp_tv + \Delta^2 v + \sumt p_{\beta}\ppp_x^{\beta}v
= -\ppp_t\www{u} - \Delta^2\www{u} - \sumt p_{\beta}\ppp_x^{\beta}\www{u}
=: F(x,t) \quad \mbox{in $\OOO \times (0,T)$}
                                                    \eqno{(4.2)}
$$
and
$$
\ppp_{\nu}^jv = 0 \quad \mbox{on}\,\, \Gamma \times (0,T),\quad  \forall j\in \{0,1,2,3\}.
                                                        \eqno{(4.3)}
$$
Next we  construct the weight function which will be used in the Carleman estimate later. 
First we construct some domain $\OOO_1$.
For $\Gamma \subset \ppp\OOO$, we choose a bounded domain $\OOO_1$ 
with smooth boundary such that
$$
\OOO \subsetneqq \OOO_1, \quad \ooo{\Gamma} = \ooo{\ppp\OOO\cap\OOO_1}, 
\quad \ppp\OOO\setminus\Gamma \subset \ppp\OOO_1.
                                                   \eqno{(4.4)}
$$
In particular, $\OOO_1\setminus\ooo{\OOO}$ contains some non-empty 
open subset. 
We note that $\OOO_1$ can be constructed as the interior of a union  
of $\ooo{\OOO}$ and the closure of a non-empty domain 
$\widehat{\OOO}$ satisfying $\widehat{\OOO} \subset \R^n 
\setminus \ooo{\OOO}$ and $\ppp\widehat{\OOO} \cap \ppp\OOO = \Gamma$.

We choose a domain $\omega$ such that 
$\ooo{\omega} \subset \OOO_1 \setminus \ooo{\OOO}$.
Then, by \cite{Ima}, we can find a function $d\in C^4(\ooo{\OOO_1})$ such that
$$
d(x)>0 \quad \mbox{in $\OOO_1$}, \quad |\nabla d(x)| > 0 \quad \mbox{on
$\ooo{\OOO_1\setminus \omega}$}, \quad d=0 \quad \mbox{on $\ppp\OOO_1$}. 
                                     \eqno{(4.5)}
$$
We recall that we choose a domain $\OOO_0 \subset \OOO$ satisfying 
$\ooo{\ppp{\OOO_0} \cap \ppp\OOO} \subset\Gamma$ and 
$\ooo{\OOO_0} \subset \OOO \cup \Gamma$.

Next we make the zero extensions of $v$ and $F$ to $\OOO_1$.  By the same 
letters, we denote the extensions:
$$
v(\cdot,t) = 
\left\{ \begin{array}{rl}
& v(\cdot,t) \quad \mbox{for}\quad (x,t)\in \OOO, \\
& 0          \quad \mbox{for}\quad (x,t)\in\OOO_1\setminus \OOO, 
\end{array}\right.
\quad F(\cdot,t) = 
\left\{ \begin{array}{rl}
& F(\cdot,t) \quad \mbox{for}\quad (x,t)\in \OOO, \\
& 0          \quad \mbox{for}\quad (x,t)\in\OOO_1\setminus \OOO. 
\end{array}\right.
$$
By (4.2) and (4.3), we can readily verify
$$
\ppp_tv + \Delta^2 v + \sumt p_{\beta}\ppp_x^{\beta}v
= F \quad \mbox{in $\OOO_1 \times (0,T)$.}    \eqno{(4.6)}  
$$

{\bf Second Step.}

Since $\OOO_0 \subset \OOO \cup \Gamma \subset \OOO_1$, we see that 
$\ooo{\OOO_0} \subset \OOO_1$.  By $d = 0$ on $\ppp\OOO_1$ in 
(4.5), we can find a small constant $\delta>0$ such that 
$$
d(x) \ge \delta, \quad x \in \ooo{\OOO_0}.     
$$
We choose small $\tau>0$ such that 
$$
\ep > \tau.                           
$$
We fix $t_0 \in (\ep, T-\ep)$ arbitrarily.
We notice that $(t_0-\tau,\, t_0+\tau)
\subset (\ep-\tau, \, T-\ep+\tau) \subset (0, T)$.

We remark that the constants $C>0$ and $\la>0$, $s_0>0$ are 
independent of $t_0$, although it is dependent on $\tau$.

Let $N>1$.  We set 
$$
h(t):= \frac{1}{ \sqrt{(t-(t_0-\tau))(t_0+\tau-t)}  },\quad
t_0-\tau < t < t_0+\tau
$$
and write $\alpha_{\tau,t_0} = \alpha$ and 
$\va_{\tau,t_0} = \va$.
Then 
\begin{align*}
& h(t) \ge \frac{1}{\tau}, \quad t_0-\tau \le t \le t_0+\tau,\\
& h(t) \le \frac{N}{\sqrt{N^2-1}}\frac{1}{\tau}, \quad
t_0-\frac{1}{N}\tau \le t \le t_0+\frac{1}{N}\tau
\end{align*}
and
$$
\alpha(x,t) = h(t)(e^{\la d(x)} - e^{2\la \dddOO})
= h(t)(1 - e^{2\la \dddOO})
$$
$$
\le \frac{1-e^{2\la \dddOO}}{\tau} =: \delta_1, \quad 
x\in \ppp\OOO_1, \, t_0-\tau \le t \le t_0+\tau.             \eqno{(4.7)}
$$
Moreover, by (4.4), we have
$$
\alpha(x,t) \ge \frac{N}{\sqrt{N^2-1}}\frac{1}{\tau}
(e^{\la \delta} - e^{2\la \dddOO}) =: \delta(N), \quad
x\in \ppp\OOO_0, \, t_0-\frac{1}{N}\tau \le t \le t_0+\frac{1}{N}\tau.
                                                 \eqno{(4.8)}
$$
Noting that $e^{\la \delta} - e^{2\la \dddOO} < 0$, we see 
that $\delta(N)$ is increasing in $N>1$.

We directly see that $\delta_1 < \delta(N)$ for each $N>1$.
Hence
$$
\delta_1 < \delta(2) < \delta(3) < \delta(4).    \eqno{(4.9)}
$$
Now we define a cut-off function $\chi \in 
C^{\infty}(\R^n\times [0,T])$ such that $0\le \chi(x,t) \le 1$ and
$$
\chi(x,t) = 
\left\{ \begin{array}{rl}
1, \quad \mbox{for}\quad (x,t)\in \{(x,t);\, \alpha(x,t) \ge \delta(3)\}, \\
0, \quad \mbox{for} \quad (x,t)\in\{(x,t);\, \alpha(x,t) \le \delta(2)\}.
\end{array}\right.
$$
We set 
$$
w = \chi v \quad \mbox{in $\OOO_1 \times (t_0-\tau,\, t_0+\tau)$}.
$$
Then
$$
\ppp_tw + \Delta^2w + \sumt p_{\beta}\ppp_x^{\beta}w 
= \chi F + G,    \eqno{(4.10)}       
$$
where 
$$
\vert G(x,t)\vert \le
\left\{ \begin{array}{rl}
& C\sums \vert \ppp_x^{\beta}v(x,t)\vert, \quad \mbox{if 
$(x,t)\in \{(x,t);\,\delta(2)<\alpha(x,t)<\delta(3)\}$}, \\  
&0, \qquad \mbox{otherwise}.
\end{array}\right.
                                             \eqno{(4.11)}
$$
In other words, 
$$
G(x,t) \ne 0 \quad \mbox{implies $\alpha(x,t) \le \delta(3)$}.
                                                         \eqno{(4.12)}
$$
By (4.7), 
$$
\alpha(x,t) \le \delta_1 \quad \forall (x,t) \in \ppp\OOO_1 \times 
(t_0-\tau, \, t_0+\tau).
$$
In view of the continuity of $\alpha$ and $\delta_1 < \delta(2)$, we can 
choose a small neighborhood $U$ of $\ppp\OOO_1$ such that 
$$
\alpha(x,t) \le \delta(2), \quad (x,t) \in U \cap \OOO_1, \,
t_0-\tau \le t \le t_0+\tau,
$$
that is, $w(x,t) = \chi(x,t)v(x,t) = 0$ for 
$(x,t) \in (U \cup \OOO_1) \times (t_0-\tau, \, t_0+\tau)$.
Therefore 
$$
w = \Delta w = 0 \quad \mbox{on $\ppp\OOO_1 \times
(t_0-\tau, \, t_0+\tau)$}.
$$
We set 
$$
J(w)(x,t) = \sumt (\vert \ppp_x^{\beta}w \vert + \vert \nabla(\Delta w)\vert
+ \vert \Delta^2w\vert + \vert \ppp_tw\vert)(x,t).
$$
Then, noting (4.11), (4.12) and $w\vert_{\omega\times (t_0-\tau, \, t_0+\tau)} = 0$, we can 
apply estimate (2.2) to equation (4.10) in $\OOO_1 \times (t_0-\tau, \, t_0+\tau)$:
\begin{align*}
& \int_{\OOO_1 \times (t_0-\tau, \, t_0+\tau)}
\frac{1}{s\va}\vert J(w)\vert^2 \weight dxdt\\
\le& C\int_{\OOO_1 \times (t_0-\tau, \, t_0+\tau)} \vert F\vert^2 \weight dxdt
+ C\int_{\{(x,t);\, \alpha(x,t) \le \delta(3)\}} 
\sums \vert \ppp_x^{\beta} v\vert^2
\weight dxdt.
\end{align*}
We set $M:= \Vert u\Vert_{L^2(0,T;H^3(\OOO))}$.
Since $\alpha(x,t) \le 0$, by (4.1) and the definition of $F$ in (4.2), 
we can estimate
$$
\int_{\OOO_1 \times (t_0-\tau, \, t_0+\tau)}
\vert F\vert^2 \weight dxdt \le CD^2
$$
and
\begin{align*}
& \int_{\{(x,t);\,\alpha(x,t) \le \delta(3)\}} 
\sums \vert \ppp_x^{\beta} v\vert^2 \weight dxdt 
\le Ce^{2s\delta(3)} \Vert v\Vert^2_{L^2(0,T;H^3(\OOO))}\\
\le & Ce^{2s\delta(3)} (\Vert u\Vert^2_{L^2(0,T;H^3(\OOO))}
+ \Vert \www{u}\Vert^2_{L^2(0,T;H^3(\OOO))})
\le Ce^{2s\delta(3)}(M^2+D^2).
\end{align*}
Hence
$$
\int_{\OOO_1 \times (t_0-\tau, \, t_0+\tau)} \frac{1}{s\va} \vert J(w)\vert^2
\weight dxdt 
\le CD^2 + Ce^{2s\delta(3)}M^2
$$
for all large $s>0$.  Shrinking the integration area 
$\OOO_1 \times (t_0-\tau, \, t_0+\tau)$ to $\OOO_0 \times 
\left( t_0 - \frac{1}{4}\tau, \, t_0 + \frac{1}{4}\tau\right)$ where 
$\alpha(x,t) \ge \delta(4)$ by (4.8), and using the fact that $\chi=1$ if 
$(x,t)\in\{(x,t);\, \alpha(x,t) \ge \delta(4)\}$, we obtain
$$
\frac{C}{s}e^{2s\delta(4)} 
\int_{\OOO_0 \times \left( t_0 - \frac{1}{4}\tau, \, t_0 + \frac{1}{4}\tau
\right)} \vert J(v)\vert^2 dxdt
\le CD^2 + Ce^{2s\delta(3)}M^2
$$
for all large $s>0$.  

Again by $u = v + \www{u}$ and (4.1), we reach 
$$
e^{2s\delta(4)} 
\int_{\OOO_0 \times \left( t_0 - \frac{1}{4}\tau, \, t_0 + \frac{1}{4}\tau
\right)} \vert J(u)\vert^2 dxdt
\le CsD^2 + Cse^{2s\delta(3)}M^2,
$$
that is,
$$
\Vert J(u)\Vert_{L^2(\OOO_0 \times \left( t_0 - \frac{1}{4}\tau, \, t_0 
+ \frac{1}{4}\tau\right)}^2 \le Ce^{Cs}D^2 + Cse^{-2s\delta_0}M^2\quad \forall s\ge s_0. \eqno{(4.12)}
$$
  
Here we set $\delta_0:= \delta(4) - \delta(3)$ and by (4.9) we note that 
$\delta_0>0$ is independent of $t_0 \in (\ep, T-\ep)$.

Since $\sup_{s>0} (s+s_0)e^{-(s+s_0)\delta_0} < \infty$, 
replacing $C>0$ by $Ce^{Cs_0}$
and changing $s$ into $s+s_0$ with $s\ge 0$ in inequality (4.12), we obtain
$$
\Vert J(u)\Vert^2_{L^2(\OOO_0\times \left(
t_0-\frac{1}{4}\tau,\, t_0+\frac{1}{4}\tau\right))} 
\le Ce^{Cs}D^2 + Ce^{-s\delta_0}M^2  
\quad \forall s\ge 0.                            \eqno{(4.13)}
$$
We minimize the right-hand side by choosing an appropriate value of 
parameter $s \ge 0$.
\\
{\bf Case 1: $M > D$.} Then we can solve 
$$
M^2e^{-s\delta_0} = e^{Cs}D^2, \quad \mbox{that is,}\quad
s = \frac{2}{C+\delta_0}\log \frac{M}{D} > 0,
$$
and, substituting into (4.13), we reach  
$$
\Vert J(u)\Vert_{L^2(\OOO_0\times \left(
t_0-\frac{1}{4}\tau,\, t_0+\frac{1}{4}\tau\right))} 
\le CM^{1-\kappa}D^{\kappa},
$$
where $\kappa = \frac{\delta_0}{C+\delta_0} \in (0,1)$.
\\
{\bf Case 2: $M \le D$.} 
Then 
$$
\Vert J(u)\Vert_{L^2(\OOO_0\times \left(
t_0-\frac{1}{4}\tau,\, t_0+\frac{1}{4}\tau\right))} 
\le C(1 + e^{Cs})^{\frac{1}{2}}D.
$$
Therefore
\begin{align*}
&\Vert u\Vert_{L^2\left( t_0-\frac{1}{4}\tau,\, t_0+\frac{1}{4}\tau;
H^2(\OOO_0)\right)}
+ \Vert \nabla(\Delta u)\Vert
_{L^2\left( t_0-\frac{1}{4}\tau,\, t_0+\frac{1}{4}\tau;
L^2(\OOO_0)\right)}\\
+ & \Vert \Delta^2u \Vert
_{L^2\left( t_0-\frac{1}{4}\tau,\, t_0+\frac{1}{4}\tau;
L^2(\OOO_0)\right)}
+ \Vert \ppp_tu \Vert
_{L^2\left( t_0-\frac{1}{4}\tau,\, t_0+\frac{1}{4}\tau;
L^2(\OOO_0)\right)}
\le C_0(M)(D^{\kappa} + D).
\end{align*}
Since the constants $C_0$ and $\kappa$ are independent of choices 
of $t_0 \in (\ep, T-\ep)$, noting that $\frac{1}{4}\tau < \tau < \ep$ and
$T-\ep < T - \tau < T-\frac{1}{4}\tau$, changing $t_0$ to divide the interval 
$(\ep, \, T-\ep)$ into a finite number of intervals in the forms
$(t_0-\tau, \, t_0+\tau)$.
Thus the proof of Theorem 2 is completed. $\square$
\section{Concluding Remarks}

{\bf 5-1.} 

For a fourth-order partial differential equation, we considered an inverse
source problem by interior data and the problem of unique continuation of solution 
from lateral Cauchy data. For both problems  stability estimates are obtained.
We can similarly argue for more general fourth-order equations and 
also inverse coefficient problems.  One related physical model equation is
$$
\ppp_ty(x,t) + \Delta^2y(x,t) - \mbox{div}\, (p(\nabla y)) = F
$$
with some function $p$, which appears for the analysis of 
epitaxial growth of nanoscale thin films (e.g., King, Stein and 
Winkler \cite{KSW}).

{\bf 5-2.}

The article relies on the Carleman estimate in Guerrero and Kassab \cite{GK}
which globally holds in the whole domain $\OOO \times (0,T)$ and is 
with the singular weight function.  In general we have two types of Carleman 
estimates: global and local.   For our fourth-order equations, we can expect 
a local one but we do not know so far.  For single partial differetial 
equations, there is a general theory providing a sufficient condition 
for a Carleman estimate (Isakov \cite{Isa}), but this fourth-order 
partial differential equation does not satisfy his sufficient condition,
in general.

Local Carleman estimates do not require 
the boundary condition on the whole $\ppp\OOO \times (0,T)$, so that the 
solutions are estimated in some subdomain, not in $\OOO\times (0,T)$.
Thus the local Carleman estimates are convenient for proving the local 
stability for inverse problems and the continuation.
However, also the global Carleman estimate can derive such local stability 
results by the argument in Section 4.
As such an earlier work, we refer to Imanuvilov and Yamamoto \cite{IY3}
for the parabolic equation, and the proof of Theorem 2 is similar to 
\cite{IY3}.

To sum up, we compare the global and the local Carleman estimates:
\begin{itemize}
\item
Global Carleman estimate (e.g., \cite{GK}, \cite{Ima})
can produce the global Lipschiz stability 
as well as local stability results.  However they do not hold in general for
hyperbolic types of equations.
\item
Local Carleman estimates are extremely difficult for deriving the global 
Lipschitz stability for parabolic types of equations.
However they hold for hyperbolic types of equations.
\end{itemize}

\section*{Acknowledgments}
The first author was partially supported by NSF grant DMS 1312900 and
Grant-in-Aid for Scientific Research (S) 
15H05740 of Japan Society for the Promotion of Science. 
The second author was supported by Grant-in-Aid for Scientific Research (S) 
15H05740 of Japan Society for the Promotion of Science and
by The National Natural Science Foundation of China 
(no. 11771270, 91730303).
This work was 
prepared with the support of the ``RUDN University Program 5-100".


\end{document}